\newif\ifMS
\MSfalse

\ifMS
\documentclass[pdflatex,sn-mathphys-num]{sn-jnl}
\else
\documentclass[letterpaper]{article}
\usepackage[margin=1.2in]{geometry}
\fi

\usepackage{mathtools}
\usepackage{amsmath,amssymb,amsfonts}%
\usepackage{amsthm}%
\usepackage{mathrsfs}%
\usepackage{textcomp}%
\usepackage{manyfoot}%
\usepackage{xparse}
\usepackage[normalem]{ulem}

\ifMS
\else
\usepackage{hyperref}
\hypersetup{
	colorlinks,
	linkcolor=blue,
	citecolor=blue,
	urlcolor=blue,
}
\fi
\usepackage[capitalize]{cleveref}

\DeclarePairedDelimiterXPP{\opnorm}[1]{}{\lVert}{\rVert}{_{\ell_2}}{#1}
\DeclarePairedDelimiterXPP{\nucnorm}[1]{}{\lVert}{\rVert}{_{*}}{#1}

\newcommand{\st}{\text{ s.t. }}
\newcommand{\Mspace}{\R^{d_1 \times d_2}}
\newcommand{\PT}{\scrP_{\scrT}}
\newcommand{\PTp}{\scrP_{\scrT^\perp}}

\makeatletter
\@ifpackageloaded{unicode-math}{%
	\newcommand{\mymathbold}{\symbf}%
}{%
	\usepackage{bm}%
	\newcommand{\mymathbold}{\bm}%
}
\makeatother

\newcommand{\scrbar}[1]{\overline{\mathcal{#1}}}
\newcommand{\scrhat}[1]{\widehat{\mathcal{#1}}}
\newcommand{\scrtl}[1]{\widetilde{\mathcal{#1}}}

\ExplSyntaxOn

\cs_new_protected:Nn \bamboo_define:nnnnnN
{
	\cs_new_protected:cpx { #3 #1 #4 } { \exp_not:N #5{#6{#2}} }
}

\int_step_inline:nnn { `A } { `Z }
{
	\bamboo_define:nnnnnN
	{ \char_generate:nn { #1 } { 11 } }
	{ \char_generate:nn { #1 } { 11 } }
	{ bl }
	{ }
	{ \mymathbold }
	\use:n
	
	\bamboo_define:nnnnnN
	{ \char_generate:nn { #1 } { 11 } }
	{ \char_generate:nn { #1 } { 11 } }
	{ scr }
	{ }
	{ \mathcal }
	\use:n
	
	\bamboo_define:nnnnnN
	{ \char_generate:nn { #1 } { 11 } }
	{ \char_generate:nn { #1 } { 11 } }
	{ }
	{ hat }
	{ \widehat }
	\use:n
	
	\bamboo_define:nnnnnN
	{ \char_generate:nn { #1 } { 11 } }
	{ \char_generate:nn { #1 } { 11 } }
	{ scr }
	{ hat }
	{ \scrhat }
	\use:n
	
	\bamboo_define:nnnnnN
	{ \char_generate:nn { #1 } { 11 } }
	{ \char_generate:nn { #1 } { 11 } }
	{ }
	{ br }
	{ \overline }
	\use:n
	
	\bamboo_define:nnnnnN
	{ \char_generate:nn { #1 } { 11 } }
	{ \char_generate:nn { #1 } { 11 } }
	{ }
	{ tl }
	{ \widetilde }
	\use:n
	
	\bamboo_define:nnnnnN
	{ \char_generate:nn { #1 } { 11 } }
	{ \char_generate:nn { #1 } { 11 } }
	{ scr }
	{ br }
	{ \scrbar }
	\use:n
	
	\bamboo_define:nnnnnN
	{ \char_generate:nn { #1 } { 11 } }
	{ \char_generate:nn { #1 } { 11 } }
	{ scr }
	{ tl }
	{ \scrtl }
	\use:n
	
	\bamboo_define:nnnnnN
	{ \char_generate:nn { #1 } { 11 } }
	{ \char_generate:nn { #1 } { 11 } }
	{ }
	{ dt }
	{ \dot}
	\use:n
}
\int_step_inline:nnn { `a } { `z }
{
	\bamboo_define:nnnnnN
	{ \char_generate:nn { #1 } { 11 } }
	{ \char_generate:nn { #1 } { 11 } }
	{ bl }
	{ }
	{ \mymathbold }
	\use:n
	
	\bamboo_define:nnnnnN
	{ \char_generate:nn { #1 } { 11 } }
	{ \char_generate:nn { #1 } { 11 } }
	{ }
	{ hat }
	{ \hat }
	\use:n
	
	\bamboo_define:nnnnnN
	{ \char_generate:nn { #1 } { 11 } }
	{ \char_generate:nn { #1 } { 11 } }
	{ }
	{ br }
	{ \bar }
	\use:n
	
	\bamboo_define:nnnnnN
	{ \char_generate:nn { #1 } { 11 } }
	{ \char_generate:nn { #1 } { 11 } }
	{ }
	{ tl }
	{ \tilde }
	\use:n                            
	
	\bamboo_define:nnnnnN
	{ \char_generate:nn { #1 } { 11 } }
	{ \char_generate:nn { #1 } { 11 } }
	{ }
	{ dt }
	{ \dot}
	\use:n
}
\clist_map_inline:nn
{
	Gamma,Delta,Theta,Lambda,Xi,Pi,Sigma,Phi,Psi,Omega
}
{
	\bamboo_define:nnnnnN
	{ #1 }
	{ #1 }
	{ bl }
	{ }
	{ \mymathbold }
	\use:c
	
	\bamboo_define:nnnnnN
	{ #1 }
	{ #1 }
	{ op }
	{ }
	{ \op }
	\use:c
	
	\bamboo_define:nnnnnN
	{ #1 }
	{ #1 }
	{ scr }
	{ }
	{ \mathcal }
	\use:c
	
	\bamboo_define:nnnnnN
	{ #1 }
	{ #1 }
	{ }
	{ hat }
	{ \widehat }
	\use:c
	
	\bamboo_define:nnnnnN
	{ #1 }
	{ #1 }
	{ }
	{ br }
	{ \overline }
	\use:c
	
	\bamboo_define:nnnnnN
	{ #1 }
	{ #1 }
	{ }
	{ tl }
	{ \widetilde }
	\use:c
	
	\bamboo_define:nnnnnN
	{ #1 }
	{ #1 }
	{ }
	{ dt }
	{ \dot}
	\use:c
	
}
\clist_map_inline:nn
{
	alpha,beta,gamma,delta,epsilon,zeta,eta,theta,iota,kappa,
	lambda,mu,nu,xi,pi,rho,sigma,tau,phi,chi,psi,omega
}
{
	\bamboo_define:nnnnnN
	{ #1 }
	{ #1 }
	{ bl }
	{ }
	{ \mymathbold }
	\use:c
	
	\bamboo_define:nnnnnN
	{ #1 }
	{ #1 }
	{ }
	{ hat }
	{ \hat }
	\use:c
	
	\bamboo_define:nnnnnN
	{ #1 }
	{ #1 }
	{ }
	{ br }
	{ \bar }
	\use:c
	
	\bamboo_define:nnnnnN
	{ #1 }
	{ #1 }
	{ }
	{ tl }
	{ \tilde }
	\use:c
	
	\bamboo_define:nnnnnN
	{ #1 }
	{ #1 }
	{ }
	{ dt }
	{ \dot}
	\use:c
}

\ExplSyntaxOff

\newcommand{\rank}{\operatorname{rank}}

\newcommand{\tr}{\operatorname{tr}}
\newcommand{\argmin}{\operatornamewithlimits{arg~min}}

\DeclarePairedDelimiter{\norm}{\lVert}{\rVert}

\DeclarePairedDelimiter{\abs}{\lvert}{\rvert}

\DeclarePairedDelimiter{\braces}{\{}{\}}
\DeclarePairedDelimiter{\parens}{(}{)}
\DeclarePairedDelimiter{\brackets}{[}{]}
\DeclarePairedDelimiterX{\ip}[2]{\langle}{\rangle}{#1,#2}

\DeclarePairedDelimiterXPP{\normsub}[2]{}{\lVert}{\rVert}{_{#2}}{#1}
\DeclarePairedDelimiterXPP{\ipsub}[3]{}{\langle}{\rangle}{_{#3}}{#1,#2}

\DeclarePairedDelimiterXPP{\ipHS}[2]{}{\langle}{\rangle}{_{\mathrm{HS}}}{#1, #2}
\DeclarePairedDelimiterXPP{\normHS}[1]{}{\lVert}{\rVert}{_{\mathrm{HS}}}{#1}

\DeclarePairedDelimiterXPP{\ipF}[2]{}{\langle}{\rangle}{_{\mathrm{F}}}{#1, #2}
\DeclarePairedDelimiterXPP{\normF}[1]{}{\lVert}{\rVert}{_{\mathrm{F}}}{#1}

\DeclarePairedDelimiterXPP{\dkl}[2]{\operatorname{D_{KL}}}{(}{)}{}{#1 \: \delimsize\Vert \: #2}

\DeclarePairedDelimiterXPP{\restr}[2]{}{{}}{\vert}{_{#2}}{#1}

\newcommand{\R}{\mathbf{R}}

\newcommand{\range}{\operatorname{range}}

\newcommand{\negqquad}{\mspace{-36mu}}

\theoremstyle{plain}%
\newtheorem{theorem}{Theorem}
\newtheorem{lemma}{Lemma}%
\newtheorem{corollary}{Corollary}%

\theoremstyle{definition}%
\newtheorem{definition}{Definition}%

\ifMS
\fi

\begin{document}

\title{Low solution rank of the matrix LASSO under RIP with consequences for rank-constrained algorithms}

\ifMS
\else
\author{Andrew D.\ McRae\thanks{A.\ McRae is with the Institute of Mathematics, EPFL, Lausanne, Switzerland. E-mail: \href{mailto:andrew.mcrae@epfl.ch}{andrew.mcrae@epfl.ch}. This work was supported by the Swiss State Secretariat for Education, Research and Innovation (SERI) under contract MB22.00027}}
\maketitle
\fi

\abstract{We show that solutions to the popular convex matrix LASSO problem (nuclear-norm--penalized linear least-squares) have low rank under similar assumptions as required by classical low-rank matrix sensing error bounds. Although the purpose of the nuclear norm penalty is to promote low solution rank, a proof has not yet (to our knowledge) been provided outside very specific circumstances. Furthermore, we show that this result has significant theoretical consequences for nonconvex rank-constrained optimization approaches. Specifically, we show that if (a) the ground truth matrix has low rank, (b) the (linear) measurement operator has the matrix restricted isometry property (RIP), and (c) the measurement error is small enough relative to the nuclear norm penalty, then the LASSO solution is unique and has rank (approximately) bounded by that of the ground truth. From this, we show (a) that a low-rank--projected proximal gradient descent algorithm will converge linearly to the unique LASSO solution from any initialization, and (b) that the nonconvex landscape of the low-rank Burer-Monteiro--factored problem formulation is benign in the sense that all second-order critical points are globally optimal and yield the unique LASSO solution.}

\section{Introduction and main result}
\label{sec:intro}
Consider the problem of estimating a low-rank matrix $M_* \in \Mspace$
from (noisy) linear measurements of the form
\begin{equation}
	\label{eq:measmodel}
	y = \scrA(M_*) + \xi,
\end{equation}
where $\scrA \colon \Mspace \to \R^n$ is a linear measurement operator,
and $\xi \in \R^n$ represents noise or other measurement error.

For such a problem,
a popular choice of estimator is a solution to the nuclear-norm--regularized least-squares optimization problem, also know as the matrix LASSO:
\begin{equation}
	\label{eq:opt_orig}
	\min_{M \in \Mspace}~\frac{1}{2} \norm{y - \scrA(M)}^2 + \lambda \nucnorm{M},
\end{equation}
where $\norm{\cdot}$ is the ordinary Euclidean $\ell_2$ norm, $\nucnorm{\cdot}$ is the nuclear norm, and $\lambda > 0$ is a hyperparameter that we can choose.
Under suitable conditions on $\scrA$ and $\xi$ and the right choice of $\lambda$,
any solution\footnote{We use the term ``solution'' in general to mean any global optimum. Although the global optimum is not necessarily unique in general for the LASSO problem, it will, in fact, be unique under the additional assumptions that we will make (our main result includes uniqueness as a byproduct).} $\Mhat$ to \eqref{eq:opt_orig} is an accurate estimate of $M_*$,
where the accuracy depends on the size of the noise $\xi$.
In fact, the matrix LASSO estimator has, up to a constant factor, minimax-optimal error under certain popular statistical models (see further discussion below).

Curiously, although the entire point of the nuclear norm penalty in \eqref{eq:opt_orig} is to promote low solution rank,
and there are, as mentioned, many theoretical guarantees in the literature that LASSO solutions are good estimators of $M_*$ in the sense of error norm,
there is little in the way of theoretical guarantees that the LASSO problem has a solution whose rank is, in fact, low (i.e., comparable to that of $M_*$).
The precise rank of a solution is important for two practical reasons:
\begin{itemize}
	\item \textbf{Storage and computation:} The space required to store a solution $\Mhat$ with rank $\rhat$ is $O(\rhat(d_1 + d_2))$.
	Similarly, the cost of multiplication by $\Mhat$ scales proportionally to $\rhat$.
	If $\rhat$ is large, we could try to save space and computational cost by taking a low-rank approximation of $\Mhat$,
	but this may introduce additional error.
	\item \textbf{Algorithmic efficiency:} The optimization problem \eqref{eq:opt_orig} is convex and has $d_1 d_2$ variables.
	Although it can be solved in polynomial time, the computation may still be impractically expensive when $d_1$ and $d_2$ are large.
	If we knew that a low-rank solution \emph{existed},
	we could consider algorithms that only search over low-rank matrices;
	this is more practical for large problem sizes.%
\end{itemize}
In this paper, we partially fill the theoretical gap for solution rank by showing that, under conditions already standard in the literature, the solution is unique and has rank similar to that of the ground truth $M_*$.
Furthermore, we show that rank-resticted optimization algorithms (which scale better than a generic full-rank algorithm)
can find this solution; see \Cref{sec:lowrank_opt} for more details.

By arriving at nonconvex low-rank algorithms via the LASSO,
we retain a key advantage of the LASSO: that it is \emph{rank-agnostic}.
We do not need to know $r^*$ precisely in order to obtain optimal error bounds.
On the other hand, the existing theory for direct (unregularized) nonconvex approaches either assumes exact knowledge of $r^*$ or suffers increased sensitivity to noise if we overestimate the rank;
see \Cref{sec:relwork_landscape} for a literature review.

Our main solution rank result will depend on three properties of our measurements in \eqref{eq:measmodel}:
\begin{itemize}
	\item \textbf{Ground truth rank:} Our result will show that solutions to \eqref{eq:opt_orig} have (approximately) the same rank $r^*$ as the ground truth $M_*$.
	
	\item \textbf{Measurement restricted isometry:}
	We need some assumption about how well $\scrA$ ``measures'' (low-rank) matrices.
	We quantify this with the popular matrix \emph{restricted isometry property} (RIP).
	We say that $\scrA$ has the $(r, \delta)$-RIP if,
	for \emph{every} matrix $M \in \Mspace$ with $\rank(M) \leq r$,
	\begin{equation}
		\label{eq:rip_def}
		(1 - \delta) \normF{M}^2 \leq \norm{\scrA(M)}^2 \leq (1 + \delta) \normF{M}^2.
	\end{equation}
	
	\item \textbf{Measurement error:} We quantify the measurement error entirely via $\opnorm{\scrA^*(\xi)}$,
	where $\opnorm{\cdot}$ is the matrix operator norm, and $\scrA^* \colon \R^n \to \Mspace$ is the adjoint of $\scrA$.
\end{itemize}
These concepts are well established in the matrix sensing literature.
It is known \cite{Candes2011b,Negahban2011,Rohde2011} that, if $r^* = \rank(M_*)$, and $\scrA$ has $(r, \delta)$-RIP for%
\footnote{We write $a \lesssim b$ (equivalently, $b \gtrsim a$) to mean $a \leq C b$ for some unspecified but universal constant $C > 0$. We write $a \approx b$ to mean $a \lesssim b$ and $a \gtrsim b$ simultaneously.}
$r \gtrsim r^*$ and sufficiently small $\delta$, then any solution $\Mhat$ to \eqref{eq:opt_orig} with regularization parameter $\lambda \approx \opnorm{\scrA^*(\xi)}$ satisfies
\[
	\normF{\Mhat - M_*} \lesssim \sqrt{r^*} \opnorm{\scrA^*(\xi)}.
\]
This error is minimax-optimal (within constants) for Gaussian noise \cite{Candes2011b,Rohde2011}.
This illustrates one key advantage of the matrix LASSO: we do not need to know the true rank $r^*$ precisely beforehand to obtain the optimal ($r^*$-dependent) error rate.
In contrast, the error of unregularized low-rank methods (see \Cref{sec:relwork_landscape}) will, in general, suffer if the search rank is chosen too large.

We can now state our main result (proved in \Cref{sec:mainproof}):
\begin{theorem}
	\label{thm:cvx_rankbd}
	Suppose we have measurements of the form \eqref{eq:measmodel}.
	Let $r^* = \rank(M_*)$,
	and suppose $\scrA$ has $(2r^*, \delta^*)$-RIP for some $\delta^* > 0$.
	Choose $\lambda > 0$, and suppose that
	\begin{equation}
		\label{eq:rankbd_cond}
		\delta^* + \frac{\opnorm{\scrA^*(\xi)}}{\lambda} \leq \frac{1}{16}.
	\end{equation}
	Then, with this choice of $\lambda$, \eqref{eq:opt_orig} has a unique solution $\Mhat$ that satisfies
	\[
	\rank(\Mhat) \leq \parens*{ 1 + 25 \parens*{\delta^* + \frac{\opnorm{\scrA^*(\xi)}}{\lambda}}^2 } r^*
	\leq \frac{11}{10} r^*.
	\]
\end{theorem}
The condition \eqref{eq:rankbd_cond} is qualitatively similar to what was assumed for prior error bounds on $\Mhat$ (see above).

Note that there is some flexibility in how we choose $M_*$ and $\xi$.
For example, if the ground truth is only \emph{approximately} low-rank,
we can decompose it into an exactly low-rank matrix plus a residual
and then set $M_*$ to be the low-rank matrix and subsume the residual (through $\scrA$) into $\xi$.
Again, see \cite{Candes2011b,Negahban2011} for examples of this type of analysis.

To illustrate our result further,
we compare it to another noisy matrix recovery problem that we understand well:
the low-rank matrix \emph{denoising} problem, where the measurement operator $\scrA$ is the identity,
and we observe $Y = M_* + \Xi$ for some noise matrix $\Xi \in \Mspace$.
The nuclear-norm--penalized problem \eqref{eq:opt_orig} then becomes
\[
\min_{M \in \Mspace}~\frac{1}{2} \normF{Y - M}^2 + \lambda \nucnorm{M}.
\]
The exact solution (which is unique due to strong convexity) can be obtained by singular value soft-thresholding $Y$:
if $Y = U_Y \Sigma_Y V_Y^T$ is the singular value decomposition (SVD) of $Y$,
then the solution is $\Mhat = U_Y (\Sigma_Y - \lambda I)_+ V_Y^T$,
where $(\Sigma_Y - \lambda I)_+$ replaces any negative diagonal elements of $\Sigma_Y - \lambda I$ by zero.
If $\opnorm{\Xi} \leq \lambda$, simple singular value perturbation arguments yield $\rank(\Mhat) \leq \rank(M_*)$.
Our result \Cref{thm:cvx_rankbd} can be interpreted as extending this (modulo constants and the RIP constant $\delta^*$) to general RIP measurement operators $\scrA$.

A caveat to our result is that the RIP assumption is quite strong.
It is well-known that it can hold (with high probability) when $\scrA$ is constructed from $O(r^*(d_1 + d_2))$ random dense sensing matrices (e.g., with i.i.d.\ Gaussian entries) \cite{Recht2010}, but such an $\scrA$ is computationally cumbersome in practice.
For accurate recovery of $M_*$ in Frobenius norm,
weaker conditions and more practical measurements suffice. For example, the paper \cite{Cai2015} gives similar recovery guarantees under an $\ell_1$ RIP (i.e., $\norm{\scrA(M)}_1 \approx \normF{M}$ for low-rank $M$).
One can construct $\scrA$ from $O(r^*(d_1 + d_2))$ random rank-1 measurements and have, with high probability, this $\ell_1$ condition but not the classical $\ell_2$ RIP that we require (again, see \cite{Cai2015}).
It is not clear how to modify our proof of \Cref{thm:cvx_rankbd} to handle such a weaker assumption.
Furthermore, the results on low-rank algorithms in \Cref{sec:lowrank_opt} depend critically on the classical RIP (as does the considerable literature on which we build---see \Cref{sec:relwork_landscape}).
Finding a result similar to \Cref{thm:cvx_rankbd} with weaker assumptions on the measurements $\scrA$ is a very interesting problem for future work.

\section{Guarantees for nonconvex rank-constrained optimization}
\label{sec:lowrank_opt}
Although the nuclear-norm--regularized problem \eqref{eq:opt_orig} is convex and can be solved in polynomial time by standard solvers,
it can be quite computationally expensive to solve directly when the matrix dimensions $d_1, d_2$ are large.

However, \Cref{thm:cvx_rankbd} guarantees, in many useful cases, that $\rhat \coloneqq \rank(\Mhat)$ is small (where, again $\Mhat$ is a global optimum of \eqref{eq:opt_orig}).
Therefore, we can consider restricted-rank versions of \eqref{eq:opt_orig}.
For any $r \geq \rhat$, $\Mhat$ will also be a global optimum of
\begin{equation}
	\label{eq:opt_lowrank}
	\min_{\substack{M \in \Mspace \\ \rank(M) \leq r}}~\frac{1}{2} \norm{y - \scrA(M)}^2 + \lambda \nucnorm{M}.
\end{equation}

The rank-restricted problem \eqref{eq:opt_lowrank} is nonconvex but lends itself well to more computationally efficient algorithms than \eqref{eq:opt_orig}.
We consider two types of algorithm: a projected proximal gradient descent algorithm and the Burer-Monteiro factorized approach.
We first present theoretical guarantees (\Cref{thm:PPGD,thm:landscape}) that are independent of \Cref{thm:cvx_rankbd},
then we combine these results with \Cref{thm:cvx_rankbd} to obtain more specific results in \Cref{cor:sensing_alg}.

\subsection{Projected spectral algorithm}
First, we consider an iterative soft/hard singular value thresholding algorithm.
Given a starting point $M_0$ and a (fixed) stepsize $\eta > 0$, the algorithm has iterates
\begin{equation}
	\label{eq:PPGD_its}
	M_{t+1} = \scrP_{r, \eta \lambda}( M_t + \eta \scrA^*(y - \scrA(M_t)) ),
\end{equation}
where, for $\alpha > 0$, $\scrP_{r, \alpha}$ is the soft/hard singular value thresholding operator given by
\begin{align*}
	\scrP_{r,\alpha}(M) &= U \begin{bmatrix*}
		(\sigma_1 - \alpha)_+ & & & & & & \\
		& \ddots & & & & \\
		& & (\sigma_r - \alpha)_+ & & & \\
		& & & 0 & & \\
		& & & & \ddots & \\
		& & & & & 0
	\end{bmatrix*} V^T, \text{ where} \\
	M &= U \begin{bmatrix*}
		\sigma_1 & & & \\
		& \sigma_2 & &\\
		& & \ddots & \\
		& & & \sigma_{\rank(M)}
	\end{bmatrix*} V^T
\end{align*}
is the SVD of $M$ with $U^T U = V^T V = I_{\rank(M)}$, $\sigma_1 \geq \sigma_2 \geq \cdots \geq \sigma_{\rank(M)} > 0$, and $(x)_+ \coloneqq \max\{x, 0\}$.
In words, $\scrP_{r,\alpha}$ soft-thresholds the singular values by $\alpha$ and then keeps only the top $r$.

Note that $\scrA^*(y - \scrA(M))$ is the negative gradient of the quadratic loss function,
so \eqref{eq:PPGD_its} can be interpreted as a projected proximal gradient descent,
where the proximal operator (singular value soft thresholding) is followed by a projection onto the space of matrices of rank at most $r$.
Without the hard low-rank projection, this would be exactly the popular iterative shrinkage-thresholding algorithm (ISTA) \cite{Beck2009} with nuclear norm regularizer.

Because, at every iteration \eqref{eq:PPGD_its},
we only need to find the top $r$ singular values/vectors of the matrix $M_t + \eta \scrA^*(y - \scrA(M_t))$,
the iterates can be computed efficiently by randomized methods (see, e.g., \cite{Oh2018}).
Because such algorithms only access the matrix via matrix products,
the low rank of $M_t$ can provide further computational savings.
On the other hand, with standard ISTA, there is no guarantee that the intermediate iterates will be low-rank.

Under matrix RIP, if a low-rank global solution to \eqref{eq:opt_orig} exists,
this algorithm (with suitable stepsize $\eta$) is guaranteed to converge to it:
\begin{theorem}
	\label{thm:PPGD}
	Let $\Mhat$ be any solution to \eqref{eq:opt_orig}, and suppose $r \geq \rhat \coloneqq \rank(\Mhat)$.
	Suppose $\scrA$ has $(2r, \delta)$-RIP for some $\delta < 1/3$.
	Then, for stepsize $\eta = \frac{1}{1+\delta}$ and any $M_0 \in \Mspace$ with $\rank(M_0) \leq r$,
	the iterates \eqref{eq:PPGD_its} satisfy
	\begin{align*}
		\normF{M_t - \Mhat}^2 \leq \parens*{ \frac{2\delta}{1 - \delta} }^t \normF{M_0 - \Mhat}^2
	\end{align*}
	for all $t \geq 0$.
\end{theorem}
Thus we have linear convergence to $\Mhat$ (which must, in fact, be the only solution to \eqref{eq:opt_orig} with rank at most $r$).
We prove this in \Cref{sec:proof_PPGD} (it is the consequence of a somewhat more general result, \Cref{lem:PPGD_gen_conv}, that may be of independent interest, though the result and its proof are straightforward adaptations of other works such as \cite[Theorem 3]{Zhang2021}).

We do have some flexibility in choosing $\eta$: if $\scrA$ is $(2r, \delta_0)$-RIP, then it is $(2r, \delta)$-RIP for any $\delta \geq \delta_0$.
Thus any $\eta$ with $\frac{3}{4} < \eta \leq \frac{1}{1 + \delta_0}$ will work,
and the convergence rate in \Cref{thm:PPGD} can be modified accordingly.

\subsection{Burer-Monteiro factorization}
Another popular algorithmic approach to \eqref{eq:opt_lowrank} is via the smooth factored problem
\begin{equation}
	\label{eq:opt_factored}
	\min_{\substack{X \in \R^{d_1 \times r} \\ Y \in \R^{d_2 \times r}}}~ \frac{1}{2} \norm{y - \scrA(X Y^T)}^2 + \frac{\lambda}{2} (\normF{X}^2 + \normF{Y}^2),
\end{equation}
which is amenable to standard local solvers such as gradient descent.
This is equivalent to \eqref{eq:opt_lowrank} because, for any matrix $M$, $\nucnorm{M} = \min_{X,Y}~\frac{\normF{X}^2 + \normF{Y}^2}{2} \st M = X Y^T$.
But again, \eqref{eq:opt_factored} is nonconvex, so, in principle, a local algorithm could get stuck in a local optimum.

In fact, this will not occur under the same conditions as \Cref{thm:PPGD}.
In fact, every \emph{second-order critical point} of \eqref{eq:opt_factored} (i.e., every point where the gradient is zero and the Hessian is positive semidefinite) yields the global optimum:
\begin{theorem}
	\label{thm:landscape}
	Under the same conditions as \Cref{thm:PPGD},
	every second-order critical point $(X, Y)$ of \eqref{eq:opt_factored} satisfies $X Y^T = \Mhat$.
\end{theorem}
This follows from \Cref{thm:PPGD} by a straightforward adaptation of \cite[Theorem 2.3(a)]{Ha2020}.
For completeness, we state and prove the adapted result as \Cref{lem:PPGDtoLandscape} in \Cref{sec:PPGDtoLandscape}.
A similar result to \Cref{thm:landscape} (with a more direct proof) is \cite[Thm.\ 1.7]{Li2019}.

\Cref{thm:landscape} guarantees that local algorithms can reach the global optimum.
For example, second-order methods such as trust-regions will reach second-order critical points (see, e.g., \cite{Cartis2012}),
and gradient descent with random initialization and a suitable step size is known to converge to a local minimum \cite{Panageas2017};
\Cref{thm:landscape} ensures that these limits will be \emph{global} minima.

One limitation of \Cref{thm:PPGD,thm:landscape} is that the maximum rank $r$ must not be so big that $\scrA$ fails to have $(2r, \delta)$-RIP.
One might expect the landscape of \eqref{eq:opt_factored} to ``improve'' monotonically as $r$ increases,
but it is not clear that this is the case. Certainly, for \emph{sufficiently large} $r$ (e.g., $r \geq \min\{d_1, d_2\}$) we can make similar guarantees, but then the restriction is trivial.
What the situation is for intermediate values of $r$ is still an open question.
Beyond our specific problem, such ``overparametrization'' is a significant theoretical issue in the literature (see \Cref{sec:relwork_landscape} for further discussion).
In practice, however, this should not be too much of a problem: because $\scrA$ is known,
we ought to know (perhaps approximately) up to what rank $\scrA$ has RIP, and there would be no benefit to setting $r$ larger than this.

\subsection{Combining the results}
We summarize the implications of the rank bound \Cref{thm:cvx_rankbd} combined with the algorithmic and landscape results \Cref{thm:PPGD,thm:landscape} as follows:
\begin{corollary}
	\label{cor:sensing_alg}
	Suppose we have measurements of the form \eqref{eq:measmodel}.
	Let $r^* = \rank(M_*)$,
	and suppose $\scrA$ has $(2r^*, \delta^*)$ RIP for some $\delta^* > 0$.
	Choose $\lambda > 0$, and suppose that
	\[
	\delta^* + \frac{\opnorm{\scrA^*(\xi)}}{\lambda} \leq \frac{1}{16}.
	\]
	Then \eqref{eq:opt_orig} has a unique solution $\Mhat$.
	Furthermore, choose optimization rank parameter
	\[
	r \geq \left\lfloor \parens*{ 1 + 25 \parens*{\delta^* + \frac{\opnorm{\scrA^*(\xi)}}{\lambda}}^2 } r^*
	\right\rfloor,
	\]
	and suppose $\scrA$ has $(2r, \delta)$-RIP for some $\delta < 1/3$. Then 
	\begin{enumerate}
		\item For stepsize $\eta = \frac{1}{1 + \delta}$ and any $M_0 \in \Mspace$ with $\rank(M_0) \leq r$,
		the iterates \eqref{eq:PPGD_its} converge linearly to $\Mhat$ with error bound
		\[
		\normF{M_t - \Mhat}^2 \leq \parens*{ \frac{2\delta}{1 - \delta} }^t \normF{M_0 - \Mhat}^2.
		\]
		\item Every second-order critical point $(X, Y)$ of the Burer-Monteiro factored problem \eqref{eq:opt_factored} satisfies $X Y^T = \Mhat$.
	\end{enumerate}
\end{corollary}

\section{Related work}
\subsection{Estimator rank bounds in matrix sensing}
Compared to our problem, it is much easier to bound estimator rank in \emph{noiseless} matrix sensing because it is possible to recover the ground truth $M_*$ exactly.
There is a plethora of results (see \cite{Davenport2016,Nguyen2019} for surveys) showing that, for various conditions on the measurements $\scrA$ (including RIP) and various choices of low-rank--inducing penalty functions (including nuclear norm), if we observe $y = \scrA(M_*)$ and consider a problem of the form
\[
\min_M~\operatorname{penalty}(M) \st \scrA(M) = y,
\]
the only global optimum is $M_*$.
Clearly, this has low rank, but this fact entirely depends on the fact that we recover $M_*$ exactly.
This can be extended to the noisy case when the error ($\xi$ in \eqref{eq:measmodel}) is sparse: that is, a large fraction of measurements remains uncorrupted. See, for example, \cite{Candes2011a,Li2017} for convex optimization estimators and \cite{Ge2017,Fattahi2020,Li2020,Ma2023c,Ding2021} for the nonconvex case.
For general non-sparse error, we cannot hope to recover $M_*$ exactly,
so this approach does not work.

For general noisy matrix sensing,
the closest result we are aware of to our \Cref{thm:cvx_rankbd} is \cite[Theorem 8]{Koltchinskii2011}. Similar results (with a nonconvex rank penalty rather than a nuclear norm) can be found in \cite{Klopp2011}.
However, the estimators in these papers are not true least-squares estimators like \eqref{eq:opt_orig} and do not yield consistency ($\Mhat \to M_*$) as the noise $\xi \to 0$.
The analysis depends on the objective function being strongly convex on $\Mspace$, which is a far stronger condition than RIP on low-rank matrices.

There are also rank-recovery results (e.g., \cite{Bunea2011,Chen2013}) for the related \emph{reduced-rank regression} problem, where our observations have the form $Y = A M_* + \Xi$; now, the linear measurement operator $M \mapsto A M$ is ordinary matrix multiplication by a known matrix $A$.
However, this model is far more structured and restrictive than the general linear model \eqref{eq:measmodel}, and the results depend on this structure.

\subsection{Analogous results in sparse recovery}
Comparable results to \Cref{thm:cvx_rankbd} have long been known for the simpler sparse recovery problem with the ordinary LASSO.
Here, given a sparse vector $x_0 \in \R^d$ and noisy observations $y = A x_0 + \xi \in \R^n$ for some $n \times d$ matrix $A$,
we attempt to estimate $x_0$ by the $\ell_1$-regularized--least-squares problem
\[
\min_{x \in \R^n}~\frac{1}{2} \norm{A x - y}^2 + \lambda \norm{x}_1.
\]
Many works in the LASSO literature show that, under certain conditions, the support (set of nonzero coordinates) of any solution is (approximately) the same (or a subset of) the ground truth support.
A representative result and many further references can be found in \cite[Chapter 11]{Hastie2015} (specifically, see Theorem 11.3).
The analysis techniques depend on the fact that the support is a \emph{discrete} quantity and therefore robust to some degree of noise.
The low-rank matrix sensing problem has no equivalent discrete structure (the closest analog is the row and column subspaces, which are continuous), so we need a different analysis technique.

\subsection{Low-rank optimization for matrix recovery}
\label{sec:relwork_landscape}
There is a vast body of literature studying nonconvex factored or low-rank--constrained problems like
\begin{equation}
	\label{eq:opt_factored_unreg}
	\min_{\substack{M \in \Mspace \\ \rank(M) \leq r}}~\frac{1}{2} \norm{y - \scrA(M)}^2 = \min_{\substack{X \in \R^{d_1 \times r} \\ Y \in \R^{d_2 \times r}}}~ \frac{1}{2} \norm{y - \scrA(X Y^T)}^2.
\end{equation}
This includes many variations (e.g., the symmetric PSD case, different loss functions, or regularization for balancing) on \eqref{eq:opt_factored_unreg} which we do not discuss here.
The representative problem \eqref{eq:opt_factored_unreg} also does not include the nuclear-norm--like penalty present in our versions \eqref{eq:opt_lowrank} and \eqref{eq:opt_factored}; the vast majority of prior work does not consider such a penalty in the rank-constrained case.
This penalty is critical for our results because it ensures a low-rank global optimum to the unrestricted problem even in the presence of noise.

Many previous results consider algorithms of the form \emph{initialization + local optimization}.
These results show that, with a sufficiently accurate initial estimate (often obtained by a spectral method), a local algorithm (such as gradient descent or alternating minimization)
will converge to a global optimum or at least to a good estimate of the ground truth.
We do not attempt to cover the full body of such results (as they have a fundamentally different character than ours), but see the survey \cite{Chi2019} as well as more recent papers such as \cite{Li2020,Yonel2020,Ding2021,Chen2023} for further references.
Of particular interest is \cite{Wang2017}, which proves error rates for Gaussian noise that are minimax optimal;
however, this assumes precise knowledge of the ground truth rank $r^*$.

More relevant to our work are results on the \emph{global optimization landscape} of problems like \eqref{eq:opt_factored_unreg}. In the following, we attempt to summarize this work. See also the survey \cite{Zhang2020}.

Many papers \cite{Ge2016,Ge2017,Li2019a,Zhang2018,Zhang2019,Sun2018,Zhang2021a,Zilber2022,Zhang2023a} directly study problems of the form \eqref{eq:opt_factored_unreg} with exact measurements ($y = \scrA(M_*)$) in contexts such as matrix sensing with RIP, phase retrieval, and matrix completion.
An interesting extension of this idea is \cite{Uschmajew2020}, which studies more general rank-constrained quadratic matrix programs and recovers results for ordinary matrix sensing as a special case.

Similar results \cite{Fattahi2020,Ma2023} can be obtained when the errors $\xi$ in \eqref{eq:measmodel} are sparse (in which case it is still feasible to recover $M_*$ exactly).
An interesting further result in this direction is \cite{Ma2023c}, which, though it does not prove a global landscape result, shows almost exact convergence to the ground truth with a small and random initialization (in addition, this paper is the only result we are aware of that does not require RIP up to rank $O(r)$ in the highly overparametrized case $r \gg \rank(M_*)$).

Other works \cite{Li2019,Zhu2018a,Ha2020,Zhu2021a,Zhang2021,Bi2021,Bi2022} study the global landscapes of more general problems of the form
\[
\min_{\substack{M \in \Mspace \\ \rank(M) \leq r}} f(M) = \min_{\substack{X \in \R^{d_1 \times r} \\ Y \in \R^{d_2 \times r}}}~f(X Y^T),
\]
where $f$ is a function whose global minimum is known to have low rank and whose Hessian has certain RIP-like properties.
The paper \cite{Li2019} is particularly relevant because it explicitly considers the nuclear norm penalty and factored equivalent that we use;
their Theorem 1.7 is very similar to our \Cref{thm:landscape}.

Finally, there are many results \cite{Bhojanapalli2016,Park2017,Zhang2018a,Ma2022,Ma2023a,Ma2023b,Ma2023} showing that, in the presence of (non-sparse) noise, we can still obtain a statistically benign landscape in the sense that every local minimum (even every second-order critical point) of \eqref{eq:opt_factored_unreg} is close to the ground truth $M_*$ (where the error depends on the amount of noise).
In the limit of zero noise, these results also yield exact recovery of the ground truth (and hence globally benign landscape).
However, none of these results, to our knowledge, guarantees optimal error rates with respect to noise as we find in, for example, \cite{Candes2011b,Negahban2011,Rohde2011} for the LASSO.

A significant caveat in nearly all these results is that they fail in the highly overparametrized case
when the factorization rank $r$ is too large for RIP to hold.
The paper \cite{Ma2023} provides examples in matrix sensing with sparse errors for which the ground truth fails to be a global minimum when the optimization rank is too large (but is the only local minimum for smaller choices of $r$).
In contrast are the papers \cite{Ma2023c} (see discussion above---this considers a particular algorithm) and \cite{Zhang2024}.
The latter shows benign landscape for all sufficiently large $r$ when the objective function is strongly convex and has Lipschitz gradient.
However, this is a very strong assumption that, for example, fails whenever the number of measurements $n < d_1 d_2$, the number of degrees of freedom.

Importantly, none of the works mentioned above (initialization + convergence or global landscape) shows that we can efficiently obtain a \emph{global optimum} to the nonconvex problem in the case of general non-sparse measurement noise;
these works contain guarantees that the resulting estimate is near the ground truth, but they leave open the possibility that there is a global optimum with strictly better objective function value.
A notable exception is \cite{Luo2022a},
which, among other things, shows that \eqref{eq:opt_factored_unreg} can be free of non-global local minima even with (general) noise.
To our knowledge, this is the only such theoretical result prior to ours.
The result has a very different character from ours;
it depends fundamentally on the (Riemannian) geometry of the space of low-rank matrices.
Hence it is only valid when we choose $r = \rank(M_*)$ exactly, and it depends on the condition number $\frac{\sigma_1(M_*)}{\sigma_r(M_*)}$.
Furthermore, it says nothing about local or global optima of \eqref{eq:opt_factored_unreg} for larger values of $r$.
Our result \Cref{thm:landscape} exploits the connection between the nonconvex problem \eqref{eq:opt_factored} and the convex problem \eqref{eq:opt_orig},
and thus we show that second-order critical points of \eqref{eq:opt_factored} are globally optimal even without any rank restriction.

\section{Proof of main rank bound}
\label{sec:mainproof}
In this section, we prove \Cref{thm:cvx_rankbd}.

\subsection{Preliminary facts}
First, we note some useful facts about nuclear norms and RIP.

Recall that if $M \in \Mspace$ is a rank-$r$ matrix with SVD $M = U \Sigma V^T$ (where $U \in \R^{d_1 \times r}$, $V \in \R^{d_2 \times r}$ with $U^T U = V^T V = I_r$, and $\Sigma$ is a diagonal $r \times r$ matrix with strictly positive diagonal entries),
then
\begin{equation}
	\label{eq:nn_subgr}
	\partial \nucnorm{M} = \{ U V^T + W : W V = 0, W^T U = 0, \opnorm{W} \leq 1 \}.
\end{equation}
We will use this several times.

Next, a well-known consequence of RIP (recall \eqref{eq:rip_def}) is the following:
\begin{lemma}
	\label{lem:ip_rip}
	Suppose $\scrA$ has $(2r, \delta)$-RIP.
	Then, for any $M_1, M_2 \in \Mspace$, each with rank $\leq r$,
	\[
	\abs*{ \ip{\scrA(M_1)}{\scrA(M_2)} - \ip{M_1}{M_2} } \leq \delta \normF{M_1} \normF{M_2}.
	\]
\end{lemma}
This is easily proved from RIP by the polarization identity.

In addition, we will want to reason about what an operator with RIP does to a matrix that is not necessarily low-rank.
For this, the following is useful:
\begin{lemma}
	\label{lem:mat_decomp}
	For any matrix $M$ and integer $r \geq 1$, we can write
	\[
	M = \sum_{k \geq 1} M_k,
	\]
	where the $M_k$'s are orthogonal in the Hilbert-Schmidt inner product,
	each $M_k$ has rank at most $r$,
	and
	\[
	\sum_{k \geq 1} \normF{M_k} \leq \normF{M} + \frac{\nucnorm{M}}{\sqrt{r}}.
	\]
\end{lemma}
Many variations of this result can be found in the sparse and low-rank--matrix recovery literature (see, e.g., \cite[proof of Theorem 3.3]{Recht2010}).
It can be proved by partitioning the singular values of $M$ into blocks of size $r$ in decreasing order;
$M_k$ is the matrix formed by the singular values and vectors corresponding to the $k$th such block.
Then, for $k \geq 2$, $\normF{M_k} \leq \frac{\nucnorm{M_{k-1}}}{\sqrt{r}}$.

With these in hand, we can begin proving the main result.

\subsection{Subgradient condition}
We will prove solution uniqueness at the end.
For now, assume $\Mhat$ is \emph{any} solution to \eqref{eq:opt_orig};
we seek to bound $\rank(\Mhat)$.
By convexity, $\Mhat$ is a solution if and only if the objective function has a zero subgradient at $\Mhat$.
Noting that $\nabla_M \frac{1}{2} \norm{y - \scrA(M)}^2 = - \scrA^* (y - \scrA(M))$,
this is equivalent to
\[
\Ehat \coloneqq \frac{1}{\lambda} \scrA^* (y - \scrA(\Mhat)) \in\partial \nucnorm{\Mhat},
\]
where $\partial \nucnorm{\Mhat}$ is the subgradient of $\nucnorm{\cdot}$ at $\Mhat$.
A clear consequence of the characterization of $\partial \nucnorm{\Mhat}$ in \eqref{eq:nn_subgr} is that
\begin{equation}
	\label{eq:rankbd_SVs}
	\rank(\Mhat) \leq \abs{\{ \ell : \sigma_\ell(\Ehat) \geq 1 \}}.
\end{equation}
Thus it suffices to bound the number of large singular values of $\Ehat = \frac{1}{\lambda} \scrA^* (y - \scrA(\Mhat))$.

\subsection{The idealized nuclear norm subgradient}
To bound \eqref{eq:rankbd_SVs}, we approximate $\Ehat$ by a matrix known to have low rank.
We will produce such a matrix by considering an ``ideal'' version of \eqref{eq:opt_orig}.
If we had that $\norm{\scrA(M)} = \normF{M}$ for all $M$ and there were no measurement error $\xi$,
\eqref{eq:opt_orig} would become
\begin{equation}
	\label{eq:opt_ideal}
	\min_{M \in \Mspace}~\frac{1}{2} \normF{M - M_*}^2 + \lambda \nucnorm{M}.
\end{equation}
As already discussed below \Cref{thm:cvx_rankbd}, the solution to this problem is well-known:
if $M_*$ has SVD $M_* = U_* \Sigma_* V_*^T$,
then the unique solution to \eqref{eq:opt_ideal} is the singular value soft-thresholded matrix $M_\lambda \coloneqq U_* \Sigma_\lambda V_*^T$,
where, if the singular values of $M_*$ are $\sigma_1, \dots, \sigma_{r^*}$,
\[
\Sigma_\lambda = \begin{bmatrix*}
	(\sigma_1 - \lambda)_+ & & & \\
	& (\sigma_2 - \lambda)_+ & & \\
	& & \ddots & \\
	& & & (\sigma_{r^*} - \lambda)_+
\end{bmatrix*},
\]
where, again, we use the shorthand $(x)_+ = \max\{ x, 0 \}$.
The counterpart to $\Ehat$ is
\begin{equation}
	E_\lambda \coloneqq \frac{1}{\lambda}(M_* - M_\lambda) = \frac{1}{\lambda} U_* \begin{bmatrix*}
		\lambda \wedge \sigma_1 & & & \\
		& \lambda \wedge \sigma_2 & & \\
		& & \ddots & \\
		& & & \lambda \wedge \sigma_{r^*}
	\end{bmatrix*} V_*^T,
\end{equation}
where $a \wedge b = \min\{a, b\}$.
We will now attempt to show that $\Ehat \approx E_\lambda$.

\subsection{A matrix error bound}
We can calculate
\begin{align}
	\lambda(\Ehat - E_\lambda)
	&= \scrA^* \scrA(M_* - \Mhat) + \scrA^*(\xi) - \lambda E_\lambda \nonumber \\
	&= \scrA^* \scrA(M_* - \Mhat - \lambda E_\lambda) + \scrA^*(\xi) + \scrA^* \scrA(\lambda E_\lambda) - \lambda E_\lambda \nonumber \\
	&= \scrA^* \scrA(M_\lambda - \Mhat) + \scrA^*(\xi) + \lambda (\scrA^* \scrA(E_\lambda) - E_\lambda). \label{eq:cert_diff_decomp}
\end{align}
The most challenging part of this last expression to bound is $\scrA^* \scrA(M_\lambda - \Mhat)$.
To do so, we will prove a bound on $\norm{\scrA(\Mhat - M_\lambda)}$.

From \eqref{eq:cert_diff_decomp}, we have
\[
\lambda \ip{ \Ehat  - E_\lambda }{M_\lambda - \Mhat} = \norm{\scrA(M_\lambda - \Mhat)}^2 + \ip{ \scrA^*(\xi) + \lambda (\scrA^* \scrA(E_\lambda) - E_\lambda) }{ M_\lambda - \Mhat },  
\]
or, rearranging,
\begin{equation}
	\label{eq:cvx_opt_eq}
	\norm{\scrA(\Mhat - M_\lambda)}^2
	= \ip{ \scrA^*(\xi) + \lambda (\scrA^* \scrA(E_\lambda) - E_\lambda) }{ \Mhat - M_\lambda } + \lambda \ip{ \Ehat  - E_\lambda }{M_\lambda - \Mhat}.
\end{equation}
Before we go further,
we write (as is common in the low-rank matrix recovery literature) $\scrT$ as the subspace of $\Mspace$ of matrices of the form $U_* A^T + B V_*^T$ for any $A \in \R^{d_2 \times r^*}$ and $B \in \R^{d_1 \times r^*}$,
and we denote by $\scrT^\perp$ its orthogonal complement, given by
\[
\scrT^\perp = \{M \in \Mspace : U_*^T M = 0, M V_* = 0\}.
\]
We denote by $\PT$ and $\PTp$ the orthogonal projections on $\Mspace$ to $\scrT$ and $\scrT^\perp$ respectively.

To upper-bound \eqref{eq:cvx_opt_eq}, the monotonicity of subgradients of convex functions implies that, for \emph{any} subgradient $E_\lambda' \in \partial \nucnorm{M_\lambda}$, we have
\[
\ip{\Ehat - E_\lambda'}{M_\lambda - \Mhat} \leq 0.
\]
Recalling the characterization of subgradients in \eqref{eq:nn_subgr},
we can choose%
\footnote{There is some subtlety here if some of $M_*$'s nonzero singular values $\sigma_1, \dots, \sigma_{r^*}$ are less than $\lambda$, but one can easily check that this choice of $E_\lambda'$ is still a valid subgradient at $M_\lambda$ in this case.}
\[
E_\lambda' = E_\lambda + W
\]
for any $W \in \scrT^\perp$ with $\opnorm{W} \leq 1$ (this implies $\ip{W}{M_\lambda} = 0$).
In particular, choose $W$ such that $\ip{W}{\Mhat} = \nucnorm{\PTp(\Mhat)}$.
Then
\begin{align*}
	\ip{\Ehat}{M_\lambda - \Mhat}
	&\leq \ip{E_\lambda + W}{M_\lambda - \Mhat} \\
	&= \ip{E_\lambda}{M_\lambda - \Mhat} - \nucnorm{\PTp(\Mhat)}.
\end{align*}
Plugging this into \eqref{eq:cvx_opt_eq}, we obtain
\begin{equation}
	\label{eq:cvx_basic_ineq}
	\norm{\scrA(\Mhat - M_\lambda)}^2
	\leq \ip{\scrA^*(\xi) + \lambda(\scrA^*\scrA(E_\lambda) - E_\lambda)}{\Mhat - M_\lambda} - \lambda \nucnorm{\PTp(\Mhat)}.
\end{equation}

For brevity, write $H \coloneqq \Mhat - M_\lambda$.
We write $H = \sum_{k\geq 1} H_k$ according to \Cref{lem:mat_decomp} with $r = r^*$.
Then, \Cref{lem:ip_rip} gives
\begin{align*}
	\ip{\scrA^* \scrA(E_\lambda) - E_\lambda}{H}
	&= \sum_{k \geq 1} \ip{\scrA(E_\lambda)}{\scrA(H_k)} - \ip{E_\lambda}{H_k} \\
	&\leq \delta^* \sum_{k \geq 1} \normF{E_\lambda} \normF{H_k} \\
	&\leq \delta^* \normF{E_\lambda} \parens*{ \normF{H} + \frac{\nucnorm{H}}{\sqrt{r^*}} } \\
	&\leq \delta^* (\sqrt{r^*} \normF{H} + \nucnorm{H}).
\end{align*}
Plugging this into \eqref{eq:cvx_basic_ineq} along with the norm duality inequality $\ip{\scrA^*(\xi)}{H} \leq \opnorm{\scrA^*(\xi)} \nucnorm{H}$ yields
\begin{equation}
	\label{eq:cvx_ineq_norms1}
	\norm{\scrA(H)}^2
	\leq \delta^* \lambda \sqrt{r^*} \normF{H} + (\delta^* \lambda + \opnorm{\scrA^*(\xi)}) \nucnorm{H} - \lambda \nucnorm{\PTp(H)}.
\end{equation}
Note that, because every matrix in $\scrT$ has rank at most $2r^*$,
\[
\nucnorm{H}
\leq \nucnorm{\PT(H)} + \nucnorm{\PTp(H)} 
\leq \sqrt{2r^*} \normF{H} + \nucnorm{\PTp(H)}.
\]
This combined with \eqref{eq:cvx_ineq_norms1} gives
\begin{equation}
	\label{eq:cvx_ineq_norms2}
	\norm{\scrA(H)}^2
	\leq ((1 + \sqrt{2}) \delta^* \lambda + \sqrt{2} \opnorm{\scrA^*(\xi)} ) \sqrt{r^*} \normF{H} - ((1- \delta^*) \lambda - \opnorm{\scrA^*(\xi)} ) \nucnorm{\PTp(H)}.
\end{equation}
The left side is nonnegative, so, if $\opnorm{\scrA^*(\xi)} < (1-\delta^*) \lambda$,
we obtain
\begin{align}
	\nucnorm{H}
	&\leq \sqrt{2r^*} \normF{H} + \nucnorm{\PTp(H)} \nonumber \\
	&\leq \parens*{ \sqrt{2r^*} + \frac{ ((1 + \sqrt{2}) \delta^* \lambda + \sqrt{2} \opnorm{\scrA^*(\xi)} ) \sqrt{r^*} }{(1- \delta^*) \lambda - \opnorm{\scrA^*(\xi)}} } \normF{H} \nonumber \\
	&= \frac{(\sqrt{2} + \delta^*) \lambda}{(1- \delta^*) \lambda - \opnorm{\scrA^*(\xi)}} \sqrt{r^*} \normF{H}. \label{eq:PTp_nn_bd}
\end{align}
We can then bound, by \Cref{lem:ip_rip},
\begin{align*}
	\abs*{ \norm{\scrA(H)}^2 - \normF{H}^2 }
	&= \abs*{ \sum_{k,\ell} \ip{\scrA(H_k)}{\scrA(H_\ell)} - \ip{H_k}{H_\ell}}  \\
	&\leq \delta^* \sum_{k, \ell} \normF{H_k} \normF{H_\ell} \\
	&\leq \delta^* \parens*{ \normF{H} + \frac{1}{\sqrt{r^*}} \nucnorm{H} }^2 \\
	&\leq \parens*{ 1 + \frac{(\sqrt{2} + \delta^*) \lambda}{(1- \delta^*) \lambda - \opnorm{\scrA^*(\xi)}} }^2 \delta^* \normF{H}^2 \\
	&\leq 8 \delta^* \normF{H}^2,
\end{align*}
where the last inequality uses the assumption that $\delta^* + \frac{\opnorm{\scrA^*(\xi)}}{\lambda} \leq 1/16$.
By the same assumption, we furthermore obtain $\norm{\scrA(H)}^2 \geq \frac{1}{2} \normF{H}^2$.
Plugging this into \eqref{eq:cvx_ineq_norms2} and dropping the nuclear norm yields
\[
\norm{\scrA(H)}^2 \leq ((1 + \sqrt{2}) \delta^* \lambda + \sqrt{2} \opnorm{\scrA^*(\xi)} ) \sqrt{2 r^*} \norm{\scrA(H)},
\]
which implies
\begin{equation}
	\label{eq:AH_bound}
	\norm{\scrA(\Mhat - M_\lambda)}
	\leq \sqrt{r^*} ((2 + \sqrt{2}) \delta^* \lambda + 2 \opnorm{\scrA^*(\xi)} ).
\end{equation}

\subsection{Final rank bounds}
We are now finally ready to bound $\Ehat - E_\lambda$.
The simplest way to exploit \eqref{eq:rankbd_SVs} would be to show that $\opnorm{\Ehat - E_\lambda} < 1$.
However, this turns out to require a strong ($r^*$-dependent) bound on $\delta$ if we only assume $\scrA$ has $(2r^*, \delta^*)$-RIP.
Instead, we will upper bound the number of large singular values that $\Ehat - E_\lambda$ can have.
Let
\[
K \coloneqq \abs{\braces{ \ell : \sigma_\ell( \Ehat - E_\lambda ) \geq 1 }}.
\]
With the convention that singular values are arranged in nonincreasing order, note that
\begin{align*}
	K \leq \sum_{\ell = 1}^{K} \sigma_\ell^2(\Ehat - E_\lambda).
\end{align*}
However, we have a circular dependence on $K$.
To get something useful, we consider two cases:
first, if $K \geq r^*$, then $\sigma_\ell(\Ehat - E_\lambda) \geq 1$ for all $\ell = 1, \dots, r^*$,
so we certainly have
\[
	r^* \leq \sum_{\ell = 1}^{r^*} \sigma_\ell^2(\Ehat - E_\lambda).
\]
On the other hand, if $K < r^*$, we have
\begin{align*}
	K \leq \sum_{\ell = 1}^{K} \sigma_\ell^2(\Ehat - E_\lambda) \leq \sum_{\ell = 1}^{r^*} \sigma_\ell^2(\Ehat - E_\lambda).
\end{align*}
We can combine these two cases into the inequality
\begin{equation}
	\label{eq:rankbd_interm}
	\min\{ r^*, K \}
	\leq \sum_{\ell = 1}^{r^*} \sigma_\ell^2(\Ehat - E_\lambda).
\end{equation}
If we can show that the right-hand side is less than $r^*$,
we will have a bound on $K$.

Now, to bound the right-hand side of \eqref{eq:rankbd_interm},
note that the function
\DeclarePairedDelimiterXPP{\rstarnorm}[1]{}{\lVert}{\rVert}{_{(r^*)}}{#1}
\[
	\Mspace \ni A \mapsto \sqrt{\sum_{\ell=1}^{r^*} \sigma_\ell^2(A)} \eqqcolon \rstarnorm{A}
\]
is a matrix norm (in fact, it is a \emph{unitarily invariant} matrix norm);
see, for example, \cite[Section 7.4.7]{Horn2013} or \cite[Section IV.2]{Bhatia1997} for more details on this type of matrix norm.

The important property of $\rstarnorm{\cdot}$ for our purposes is that it has the following variational form:
\[
	\rstarnorm{A} = \max_{\substack{B \in \Mspace \\ \normF{B} \leq 1 \\ \rank(B) \leq r^*}}~\ip{A}{B}.
\]
Indeed, for $B$ of rank at most $r^*$, $A^T B$ has rank at most $r^*$, so
\begin{align*}
	\ip{A}{B}
	&= \tr(A^T B) \\
	&\leq \sum_{\ell=1}^{r^*} \sigma_\ell(A^T B) \\
	&\leq \sum_{\ell=1}^{r^*} \sigma_\ell(A) \sigma_\ell(B) \\
	&\leq \rstarnorm{A} \normF{B},
\end{align*}
where the middle inequality is a standard singular value inequality for matrix products,
and the last inequality is Cauchy-Schwartz on the singular values.
Equality is achieved for $B = \frac{1}{\rstarnorm{A}} A_{(r^*)}$, where $A_{(r^*)}$ is any (not necessarily unique) rank-$r^*$ projection of $A$.

Returning to our problem, we obtain, from \eqref{eq:rankbd_interm} and the above considerations,
\[
	\min\{ r^*, K \} \leq \max_{\substack{B \in \Mspace \\ \normF{B} \leq 1 \\ \rank(B) \leq r^*}}~\ip{\Ehat - E_\lambda}{B}^2.
\]

To bound this, note that for any $B \in \Mspace$ with $\normF{B} \leq 1$ and $\rank(B) \leq r^*$,
we have, by \eqref{eq:cert_diff_decomp},
\begin{align*}
	\abs{\ip{\Ehat - E_\lambda}{B}}
	\leq \frac{1}{\lambda} \abs{\ip{\scrA^* \scrA(M_\lambda - \Mhat)}{B}}
	+ \abs{\ip{\scrA^* \scrA(E_\lambda) - E_\lambda}{B}}
	+ \frac{1}{\lambda} \abs{\ip{\scrA^*(\xi)}{B}}.
\end{align*}
We can bound each of the three terms on the right-hand side:
\begin{itemize}
	\item The simplest is
	\[
	\abs{\ip{\scrA^*(\xi)}{B}} \leq \opnorm{\scrA^*(\xi)} \nucnorm{B} \leq \sqrt{r^*} \opnorm{\scrA^*(\xi)}.
	\]
	\item Next, by \Cref{lem:ip_rip},
	\begin{align*}
		\abs{\ip{\scrA^* \scrA(E_\lambda) - E_\lambda}{B}}
		&= \abs{\ip{\scrA(E_\lambda)}{\scrA(B)} - \ip{E_\lambda}{B} } \\
		&\leq \delta^* \normF{E_\lambda} \normF{B} \\
		&\leq \sqrt{r^*} \delta^*.
	\end{align*}
	\item Finally, using \eqref{eq:AH_bound} and the original RIP inequality \eqref{eq:rip_def},
	\begin{align*}
		\abs{\ip{\scrA^* \scrA(M_\lambda - \Mhat)}{B}}
		&\leq \norm{\scrA(M_\lambda - \Mhat)} \norm{\scrA(B)} \\
		&\leq \sqrt{r^*} ((2 + \sqrt{2}) \delta^* \lambda + 2 \opnorm{\scrA^*(\xi)})  \cdot \sqrt{1 + \delta^*} \\
		&\leq 4 \sqrt{r^*} (\delta^* \lambda + \opnorm{\scrA^*(\xi)}).
	\end{align*}
\end{itemize}
Putting these inequalities together, we obtain
\[
\min\{ r^*, K \} \leq 25 \parens*{ \delta^* + \frac{\opnorm{\scrA^*(\xi)}}{\lambda} }^2 r^*.
\]
The assumption $\delta^* + \frac{\opnorm{\scrA^*(\xi)}}{\lambda} \leq \frac{1}{16}$ ensures that the right-hand side is strictly less than $r^*$,
so the left-hand side is equal to $K$.
Because $\rank(E_\lambda) \leq r^*$,
we obtain, by \eqref{eq:rankbd_SVs},
\begin{align*}
	\rank(\Mhat)
	&\leq \abs{\{ \ell : \sigma_\ell(\Ehat) \geq 1 \}} \\
	&\leq \rank(E_\lambda) + K \\
	&\leq \brackets*{ 1 + 25 \parens*{ \delta^* + \frac{\opnorm{\scrA^*(\xi)}}{\lambda} }^2 } r^*.
\end{align*}
One can easily verify that this last bound is $\leq \frac{11}{10} r^*$ under the assumption $\delta^* + \frac{\opnorm{\scrA^*(\xi)}}{\lambda} \leq \frac{1}{16}$.
Except for solution uniqueness, this completes the proof of \Cref{thm:cvx_rankbd}.

\subsection{Uniqueness}
For any integer $p \geq 1$, we denote by $\delta_p$ the smallest possible such constant such that $\scrA$ has $(p, \delta_p)$-RIP.
Set $r = \lfloor 11 r^*/ 10 \rfloor < \frac{3}{2} r^*$.
Then $2r < 3 r^*$, and
\[
\delta_{2r} \leq \delta_{3 r^*} \leq 3 \delta_{2 r^*} \leq 3 \delta^* \leq \frac{1}{5}
\]
(the second inequality is a standard property of RIP constants---see, e.g., \cite[Proposition 6.6]{Foucart2013} for the analogous result in the sparse vector case).

Let $\Mhat, \Mhat'$ be two global optima to \eqref{eq:opt_orig}.
We have shown that $\Mhat, \Mhat'$ both have rank at most $r$.
Set, as before,
\[
\Ehat \coloneqq \frac{1}{\lambda} \scrA^*(y - \scrA(\Mhat)) \in \partial \nucnorm{\Mhat}.
\]
The fact that $\Mhat$ and $\Mhat'$ are both optimal implies
\begin{align*}
	0 &= \frac{1}{2} \norm{y - \scrA(\Mhat')}^2 + \lambda \nucnorm{\Mhat'} - \frac{1}{2} \norm{y - \scrA(\Mhat)}^2 - \lambda \nucnorm{\Mhat} \\
	&= \frac{1}{2} \norm{\scrA(\Mhat - \Mhat')}^2 + \ip{\scrA^*(\scrA(\Mhat) - y)}{\Mhat' - \Mhat} + \lambda( \nucnorm{\Mhat'} - \nucnorm{\Mhat} ) \\
	&\geq \frac{1 - \delta_{2r}}{2} \normF{\Mhat - \Mhat'}^2 + \ip{\underbrace{\scrA^*(\scrA(\Mhat) - y)}_{= - \lambda \Ehat}}{\Mhat' - \Mhat} + \lambda \ip{\Ehat}{\Mhat' - \Mhat} \\
	&= \frac{1 - \delta_{2r}}{2} \normF{\Mhat - \Mhat'}^2 \\
	&\geq \frac{2}{5} \normF{\Mhat - \Mhat'}^2.
\end{align*}
The first inequality uses $2r$-RIP and the convexity of the nuclear norm.
Clearly, this implies $\Mhat = \Mhat'$.
Thus the solution $\Mhat$ is unique.

\section{Additional proofs}
In this section, we state and prove two results we need on composite optimization that,
though they are straightforward extensions of existing results, have not, to our knowledge, appeared as such in the literature (thus they may be of independent interest).
\subsection{Convergence of constrained proximal gradient descent}
\label{sec:proof_PPGD}
\Cref{thm:PPGD} is an immediate consequence of the following more general result on constrained proximal gradient descent.
Let $\scrC \subseteq \R^d$ be a closed set,
and consider the constrained composite minimization problem
\begin{equation}
	\min_{x \in \scrC}~f(x) + h(x),
\end{equation}
where $f$ is differentiable,
and $h$ is a convex and bounded-below function such that, for any $\eta > 0$ and $x \in \R^d$, we can evaluate the  constrained proximal map
\[
\scrP_\eta(x) \coloneqq \argmin_{y \in \scrC}~\frac{1}{2\eta}\norm{x - y}^2 + h(y).
\]
If the minimum is not unique, we can choose one arbitrarily.

Given $x_0 \in \scrC$ and stepsize $\eta > 0$, we consider the iterative constrained proximal gradient descent algorithm with iterates
\begin{equation}
	\label{eq:PPGD_gen_its}
	x_{t+1} = \scrP_\eta( x_t - \eta \nabla f(x_t) ) = \argmin_{y \in \scrC}~\ip{\nabla f(x_t)}{y} + \frac{1}{2\eta}\norm{x_t - y}^2 + h(y).
\end{equation}
To analyze the performance of this algorithm, we will need to make a regularity assumption on $f$.
\begin{definition}
	$f$ is $(\scrC, \delta)$-regular if, for all $x, y \in \scrC$,
	\[
	\abs*{ f(y) - f(x) - \ip{\nabla f(x)}{y - x} - \frac{1}{2}\norm{y - x}^2} \leq \frac{\delta}{2} \norm{y - x}.
	\]
\end{definition}

We now have what we need to state our general convergence result:
\begin{lemma}
	\label{lem:PPGD_gen_conv}
	Suppose $f$ is $(\scrC, \delta)$-regular for some $\delta < 1/3$,
	and suppose there is a minimizer $x^*$ of the unconstrained problem
	\[
	\min_{x \in \R^d}~f(x) + h(x)
	\]
	such that $x^* \in \scrC$.
	Then, choosing stepsize $\eta = \frac{1}{1 + \delta}$, the iterates \eqref{eq:PPGD_gen_its} satisfy
	\begin{equation}
		\label{eq:PPGD_gen_conv}
		\norm{x_{t+1} - x^*}^2 \leq \frac{2 \delta}{1 - \delta} \norm{x_{t} - x^*}^2
	\end{equation}
	for all $t \geq 0$, and therefore $\{ x_t \}_{t \geq 0}$ converges linearly to $x^*$.
\end{lemma}
\Cref{thm:PPGD} is a special case of this, identifying $\R^d$ with $\Mspace$, taking $\scrC \subset \Mspace$ to be the set of matrices of rank at most $r$,
and setting $f(M) = \frac{1}{2}\norm{y - \scrA(M)}^2$
and $h(M) = \lambda \nucnorm{M}$.
That $f$ is $(\scrC, \delta)$-regular follows immediately from $(2r, \delta)$-RIP.

In the remainder of this section, we prove \Cref{lem:PPGD_gen_conv}.
The proof is a modification and generalization of that of \cite[Theorem 3]{Zhang2021},
with many ideas (and much of our notation) taken from \cite{Taylor2018}.

Let $\partial h(x)$ denote the subgradient set of the convex function $h$ at $x$.
Because $x^*$ minimizes $x \mapsto f(x) + h(x)$ over all $\R^d$,
we have
\begin{equation}
	\label{eq:subgrad_opt}
	s^* \coloneqq - \nabla f(x^*) \in \partial h(x^*).
\end{equation}
Now, fix an integer $t \geq 0$.
To prove \eqref{eq:PPGD_gen_conv}, we will combine \eqref{eq:subgrad_opt} with the following inequalities:
\begin{align}
	f(x_{t+1}) - f(x_t) &\leq \ip{\nabla f(x_t)}{x_{t+1} - x_t} + \frac{1+\delta}{2} \norm{x_{t+1} - x_t}^2 \label{eq:reg_ub} \\
	f(x^*) - f(x_t) &\geq \ip{\nabla f(x_t)}{x^* - x_t} + \frac{1-\delta}{2} \norm{x^* - x_t}^2 \label{eq:reg_lb} \\
	f(x_{t+1}) - f(x^*) &\geq \ip{\nabla f(x^*)}{x_{t+1} - x^*} + \frac{1-\delta}{2} \norm{x^* - x_{t+1}}^2 \label{eq:reg_lb2} \\
	\begin{split}
	&\negqquad \negqquad \ip{\nabla f(x_t)}{x_{t+1}} + \frac{1 + \delta}{2} \norm{x_{t+1} - x_t}^2 + h(x_{t+1}) \\
	&\leq \ip{\nabla f(x_t)}{x^*} + \frac{1 + \delta}{2} \norm{x^* - x_t}^2 +h(x^*) \label{eq:prox_ineq} \end{split} \\
	h(x_{t+1}) &\geq h(x^*) + \ip{s^*}{x_{t+1} - x^*}. \label{eq:ineq_subgr_convex}
\end{align}
Inequalities \eqref{eq:reg_ub}, \eqref{eq:reg_lb}, and \eqref{eq:reg_lb2} follow immediately from $f$ being $(\scrC, \delta)$-regular.
Using $\eta = \frac{1}{1+\delta}$,
\eqref{eq:prox_ineq} follows from the optimality of $x_{t+1}$ in \eqref{eq:PPGD_gen_its} and the fact that $x^* \in \scrC$.
\eqref{eq:ineq_subgr_convex} is due to the convexity of $h$.

Write $F(x) = f(x) + h(x)$ as the composite objective function.
We calculate
\begin{align*}
	F(x_{t+1}) - F(x_t)
	&= f(x_{t+1}) - f(x_t) + h(x_{t+1}) - h(x_t) \\
	&\overset{\mathclap{\eqref{eq:reg_ub}}}{\leq} \ip{\nabla f(x_t)}{x_{t+1} - x_t} + \frac{1+\delta}{2} \norm{x_{t+1} - x_t}^2 + h(x_{t+1}) - h(x_t) \\
	&\overset{\mathclap{\eqref{eq:prox_ineq}}}{\leq} \ip{\nabla f(x_t)}{x^* - x_t} + \frac{1+\delta}{2} \norm{x^* - x_t}^2 + h(x^*) - h(x_t) \\
	&\overset{\mathclap{\eqref{eq:reg_lb}}}{\leq} f(x^*) - f(x_t) - \frac{1-\delta}{2} \norm{x^* - x_t}^2 + \frac{1+\delta}{2} \norm{x^* - x_t}^2 + h(x^*) - h(x_t) \\
	&= F(x^*) - F(x_t) + \delta \norm{x^* - x_t}^2.
\end{align*}
Therefore,
\begin{equation}
	\label{eq:interm1}
	F(x_{t+1}) - F(x^*) \leq \delta \norm{x^* - x_t}^2.
\end{equation}
However, we also have
\begin{align*}
	F(x_{t+1}) - F(x^*)
	&\overset{\mathclap{\eqref{eq:reg_lb2}}}{\geq} \ip{\nabla f(x^*)}{x_{t+1} - x^*} + \frac{1-\delta}{2} \norm{x^* - x_{t+1}}^2 + h(x_{t+1}) - h(x^*) \\
	&\overset{\mathclap{\eqref{eq:subgrad_opt}}}{=} \frac{1-\delta}{2} \norm{x^* - x_{t+1}}^2 - \ip{s^*}{x_{t+1} - x^*} + h(x_{t+1}) - h(x^*) \\
	&\overset{\mathclap{\eqref{eq:ineq_subgr_convex}}}{\geq} \frac{1-\delta}{2} \norm{x^* - x_{t+1}}^2.
\end{align*}
Combining this with \eqref{eq:interm1} yields the result \eqref{eq:PPGD_gen_conv}.

\subsection{SOCPs of factored problem are stationary for PPGD}
\label{sec:PPGDtoLandscape}
In this section, we state and prove the below lemma,
which is a straightforward extension of \cite[Theorem 2.3(a)]{Ha2020} to the case of nuclear norm regularization.
A slightly weaker version of this result appeared in the blog post \cite{Boumal2024},
which was itself in part inspired by an earlier version of our paper.
 
For twice-differentiable $f$, we denote by $\nabla^2 f(x)[\xdt, \xdt]$ the Hessian quadratic form of $f$ at $x$ applied to the direction $\xdt$.

\begin{lemma}
	\label{lem:PPGDtoLandscape}
	Let $f \colon \Mspace \to \R$ be twice continuously differentiable.
	Fix a positive integer $r \leq \min\{d_1,d_2\}$.
	For $\lambda > 0$, consider the	factored low-rank nuclear-norm--regularized optimization problem
	\begin{equation}
		\label{eq:opt_factored_gen}
		\min_{\substack{X \in \R^{d_1 \times r} \\ Y \in \R^{d_2 \times r}}}~ f(X Y^T) + \frac{\lambda}{2} (\normF{X}^2 + \normF{Y}^2).
	\end{equation}
	Suppose, furthermore, that $f$ has bounded Hessian in the following sense:
	there is some $L > 0$ such that, for all $A \in \Mspace$ of rank at most $r$ and $\Adt \in \Mspace$ of rank at most $2$,
	\[
		\nabla^2 f(A)[\Adt, \Adt] \leq L \normF{A}^2.
	\]
	Then, if $0 < \eta \leq \frac{1}{L}$, every second-order critical point $(X,Y)$ of the factored problem \eqref{eq:opt_factored_gen} yields a fixed point of the rank-$r$ PPGD iteration map  with stepsize $\eta$ in the sense that $M \coloneqq X Y^T$ satisfies
	\begin{equation}
		\label{eq:PPGD_fixedpoint}
		M \in \argmin_{\substack{M' \in \Mspace \\ \rank(M') \leq r}}~\ip{\nabla f(M)}{M'} + \frac{1}{2\eta}\normF{M' - M}^2 + \lambda \nucnorm{M'}.
	\end{equation}
\end{lemma}
The Hessian assumption could easily be relaxed to a locally Lipschitz gradient assumption like in \cite{Ha2020}, but we omit this generalization for simplicity and because we do not need it.

With this lemma, \Cref{thm:landscape} follows immediately from \Cref{thm:PPGD}.
Indeed, $(2r, \delta)$-RIP implies the Hessian condition with $L = 1 + \delta$,
and the stepsize $\eta$ used in \Cref{thm:PPGD} is precisely $\eta = \frac{1}{1 + \delta} = \frac{1}{L}$.

\begin{proof}[Proof of \Cref{lem:PPGDtoLandscape}]
	Let $(X, Y)$ be a second-order critical point of \eqref{eq:opt_factored_gen}, and set $M = X Y^T$.
	By standard calculations (see, e.g., \cite{Li2019} or, for the part corresponding to $f$, \cite{Ha2020}),
	this corresponds to the following properties.
	First, the gradient of the composite objective function in $(X, Y)$ is zero, that is,
	\begin{align*}
		\nabla f(M) Y &= - \lambda X \text{ and} \\
		\nabla f(M)^T X &= - \lambda Y.
	\end{align*}
	Second, the Hessian is positive semidefinite, which means,
	for all $\Xdt \in \R^{d_1 \times r}, \Ydt \in \R^{d_2 \times r}$,
	\[
		\ip{\nabla f(M)}{\Xdt \Ydt^T} + \frac{1}{2} \nabla^2 f(M)[ X \Ydt^T + \Xdt Y^T, X \Ydt^T + \Xdt Y^T ] + \frac{\lambda}{2} (\normF{\Xdt}^2 + \normF{\Ydt}^2) \geq 0.
	\]
	
	As $\lambda > 0$, the zero-gradient conditions imply that the pair $(X, Y)$ is \emph{balanced} (see, e.g., \cite{Li2019} for more details).
	Let $r_M = \rank(M)$,
	and let $M = U \Sigma V^T$ be a compact singular value decomposition of $M$ with $U^T U = V^T V = I_{r_M}$ and diagonal $\Sigma \in \R^{r_M \times r_M}$ with strictly positive diagonal entries.
	Then, for some $r \times r_M$ matrix $R$ satisfying $R^T R = I_{r_M}$,
	we have $X = U \Sigma^{1/2} R^T$ and $Y = V \Sigma^{1/2} R^T$;
	these are valid singular value decompositions of $X$ and $Y$,
	and we see that $X$ and $Y$ have the same singular values and right singular vectors.
	Furthermore, multiplying the zero-gradient conditions on the right by $R \Sigma^{-1/2}$, we obtain
	\begin{align*}
		\nabla f(M) V &= - \lambda U \text{ and} \\
		\nabla f(M)^T U &= - \lambda V.
	\end{align*}
	We can then write
	\[
		\nabla f(M) = - \lambda U V^T + W^\perp,
	\]
	where $W^\perp$ is a matrix satisfying $W^\perp V = 0$ and $U^T W^\perp = 0$.
	
	This is almost all we need to show \eqref{eq:PPGD_fixedpoint} (cf.\ \cite[eq.\ (2.4)]{Ha2020}).
	The only additional optimality condition we need is
	\begin{equation}
		\label{eq:Wperp_bd}
		\opnorm{W^\perp} \leq \lambda + \frac{1}{\eta} \sigma_r(M).
	\end{equation}
	To show this, we apply the Hessian condition.
	Due to the fact that $(X,Y)$ is balanced as discussed above,
	there exists a unit-norm $v \in \R^r$ such that $\norm{Xv}^2 = \norm{Yv}^2 = \sigma_r(M)$.
	For arbitrary unit-norm $x \in \range(U)^\perp$, $y \in \range(V)^\perp$,
	taking $\Xdt = x v^T$, $\Ydt = y v^T$ in the Hessian condition gives
	\begin{align*}
		0 &\leq \ip{\nabla f(M)}{x y^T} + \frac{1}{2} \nabla^2 f(M)[ X v y^T + x (Yv)^T, X v y^T + x (Yv)^T ] + \lambda \\
		&\leq \ip{\nabla f(M)}{x y^T} + \frac{L}{2} \normF{X v y^T + x (Yv)^T}^2 + \lambda \\
		&= \ip{W^\perp}{x y^T} + L \sigma_r(M) + \lambda.
	\end{align*}
	As this holds for any $x \in \range(U)^\perp$, $y \in \range(V)^\perp$,
	we must have $\opnorm{W^\perp} \leq \lambda + L \sigma_r(M)$.
	The condition $\eta \leq \frac{1}{L}$ ensures the optimality condition \eqref{eq:Wperp_bd} holds,
	so the result \eqref{eq:PPGD_fixedpoint} follows.
	
\end{proof}

\ifMS
\else
\bibliographystyle{ieeetr}
\fi

\bibliography{refs}

\end{document}